\documentclass[12pt]{amsart}

\usepackage{amsmath,amssymb,amscd,amsthm,epsf,graphicx}

\numberwithin{equation}{section}
\theoremstyle{plain}
\newtheorem{theorem}[equation]{Theorem}

\newtheorem{lemma}[equation]{Lemma}
\newtheorem{corollary}[equation]{Corollary}

\theoremstyle{remark}

\theoremstyle{definition}
\newtheorem{definition}[equation]{Definition}

\newcommand{\ra}{\rightarrow}

\newcommand{\C}{\mathbb{C}}

\renewcommand{\H}{\mathbb H}

\newcommand{\R}{\mathbb R}

\newcommand{\acts}{\curvearrowright}

\newcommand{\qi}{\operatorname{QI}}

\newcommand{\al}{\alpha}

\def\Ga{\Gamma}

\def\ra{\rightarrow}

\def\si{\sigma}

\def\be{\beta}

\def\geo{\partial}

\def\defeq{:=}

\def\XXint#1#2#3{{\setbox0=\hbox{$#1{#2#3}{\int}$}
     \vcenter{\hbox{$#2#3$}}\kern-.5\wd0}}

\begin{document}

\title{Induced quasi-actions: a remark}

\author{Bruce Kleiner}
\address{Bruce Kleiner\,: Mathematics Department\\
         Yale University\\
             New Haven, CT 06520}
\email{bruce.kleiner@yale.edu}
\thanks{The first author was partially supported by NSF Grant
DMS 0701515. }
\author{Bernhard Leeb}
\address{Bernhard Leeb:
Math. Anst., Univ. M\"unchen\\
Theresienstr. 39\\
D-80333 M\"unchen}
\email{ b.l@lmu.de}

\date{July 30, 2007}
\maketitle

\section{Introduction}

In this note we observe that the notion of an induced representation
has an analog for quasi-actions, and give some applications.

We will use the definitions and notation from \cite{gpsqiss}.

\subsection{Induced quasi-actions and their properties}
Let $G$ be a group and $\{X_i\}_{i\in I}$ be a finite collection
of unbounded metric spaces. 
\begin{definition}
A quasi-action  $G\stackrel{\rho}{\acts} \prod_i\;X_i$ 
{\em preserves the product structure} if  each $g\in G$ acts by a product
of quasi-isometries, up to uniformly bounded error.
Note that 
we allow the quasi-isometries $\rho(g)$ to permute the factors,
i.e.\ $\rho(g)$ is uniformly close to a map of the form 
$(x_i)\mapsto\bigl(\phi_{\si^{-1}(i)}(x_{\si^{-1}(i)})\bigr)$
with a permutation $\si$ of $I$ and 
quasi-isometries $\phi_i:X_i\mapsto X_{\si(i)}$.
\end{definition}

Associated to every  quasi-action 
$G\stackrel{\rho}{\acts} \prod_i\;X_i$ preserving product structure
is the  action 
$G\stackrel{\rho_I}{\acts} I$ corresponding to 
the induced permutation of the factors; this is well-defined because the $X_i$'s
are unbounded metric spaces.
For each $i\in I$, 
the stabilizer $G_i$ of $i$ with respect to $\rho_I$ 
has a quasi-action $G_i\acts X_i$ by restriction of $\rho$.
It is well-defined up to equivalence 
in the sense of \cite[Definition 2.3]{gpsqiss}.

If the permutation action $\rho_I$ is {\em transitive},
all factors $X_i$ are quasi-isometric to each other,
and the restricted quasi-actions $G_i\acts X_i$ 
are quasi-conjugate
(when identifying different stabilizers $G_i$ by inner automorphisms of $G$).
The main result of this note is
that in this case any of the quasi-actions $G_i\acts X_i$ 
determines $\rho$ up to quasi-conjugacy,  and moreover
any quasi-conjugacy class may arise as 
a restricted action.

\begin{theorem}
\label{thminducedproperties}
Let $G$ be a group,  $H$ be a finite index subgroup,  and
 $H\stackrel{\al}{\acts} X$ be a quasi-action of $H$ on an unbounded 
metric space $X$.   
Then there exists a quasi-action 
$G\stackrel{\be}{\acts} \prod_{i\in G/H}\; X_i$ 
preserving product structure,
where 
\begin{enumerate}
\item Each factor $X_i$ is quasi-isometric to $X$.

\item The associated action $G\stackrel{\be_{G/H}}{\acts}\;G/H$ is 
 the natural action by left multiplication. 

\item The restriction of $\be$ to a quasi-action of $H$ on $X_H$
is quasi-conjugate to $H\stackrel{\al}{\acts}\; X$.
\end{enumerate}
Furthermore, there is a unique such quasi-action $\be$
preserving the product structure, 
up to  quasi-conjugacy by a product quasi-isometry.
Finally, if $\al$ is an isometric action, 
then the $X_i$ may be taken isometric to $X$ 
and $\be$ may be taken to be an isometric action.
\end{theorem}

\begin{definition}
\label{definducedquasiaction}
Let $G$, $H$ and $H\acts X$ be as in Theorem \ref{thminducedproperties}.
The quasi-action $\be$ is called the {\em quasi-action  induced by $H\acts X$}.
\end{definition}

As a byproduct of the main construction, we get the following:
\begin{corollary}
If $G\stackrel{\rho}{\acts} X$ is an $(L,A)$-quasi-action on an arbitrary metric space $X$,
then $\rho$ is $(L,3A)$-quasi-conjugate to a canonically defined
isometric action $G\acts X'$.
\end{corollary}

\subsection{Applications}
 The implication of Theorem \ref{thminducedproperties} 
is that in order to
quasi-conjugate a quasi-action on a product to an isometric action, 
it suffices
to quasi-conjugate the factor quasi-actions to isometric actions.
We begin with a special case:

\begin{theorem}
\label{thmcobounded}
Let $G\stackrel{\rho}{\acts} X$ be a cobounded quasi-action on $X=\prod_i \;X_i$,
where each $X_i$ is either an irreducible symmetric space of noncompact 
type, or a thick irreducible Euclidean building of rank at least two,
with cocompact Weyl group.  
Then $\rho$ is quasi-conjugate to an isometric action on $X$, after
suitable rescaling of the metrics on the factors $X_i$.
\end{theorem}

{\em Remarks}
 
\begin{itemize}
\item Theorem \ref{thmcobounded}
 was stated incorrectly
as Corollary 4.5 in \cite{gpsqiss}.  The proof
given there was was only valid for quasi-actions which do not permute the factors. 
\item Rescaling of the factors
is necessary, in general: if one takes the product of two copies of $\H^2$ where the
factors are scaled to have different curvature, then a quasi-action which 
permutes the factors will not be quasi-conjugate to an isometric action. 
\end{itemize}

We now consider a more general situation.
Let $G\stackrel{\al}{\acts}\;\prod_{i\in I}\;X_i$ be a 
quasi-action, 
where each $X_i$ is one of the following four types of spaces: 
\begin{enumerate}
\item An irreducible symmetric space of noncompact type.
\item A thick irreducible Euclidean building of rank/dimension $\geq2$,
with cocompact Weyl group. 
\item A bounded valence bushy tree in the sense of \cite{moshersageevwhyte}. 
We recall that a tree is {\em bushy} if each of its points lies within uniformly
bounded distance from a vertex having at least three unbounded
complementary components.
\item A quasi-isometrically rigid Gromov hyperbolic space which is 
of coarse type I in the sense of \cite[sec.\ 3]{kkl} 
(see the remarks below).  
A space is {\em quasi-isometrically rigid} 
if every $(L,A)$-quasi-isometry is at
distance at most $D=D(L,A)$ from a unique isometry.
\end{enumerate}
By \cite[Theorem B]{kkl}, the quasi-action preserves product structure,
and hence we have an induced permutation action $G\acts I$.
Let $J\subset I$ be the set of indices $i\in I$ such that $X_i$
is either a real 
hyperbolic space $\H^k$ for some $k\geq 4$,  a complex
hyperbolic space $\C\H^l$ for some $l\geq 2$, 
or a bounded valence bushy tree.  Generalizing
Theorem \ref{thmcobounded} we obtain:
\begin{theorem}
\label{thmquasiaction}
If the quasi-action $G_j\acts X_j$ is cobounded for each $j\in J$, 
then 
$\al$ is quasi-conjugate by a product quasi-isometry to an isometric action 
$G\stackrel{\al'}{\acts}\;\prod_{i\in I}\;X_i'$,
where for every $i$,  $X_i'$ is quasi-isometric to $X_i$,
and precisely one of the following holds:
\begin{enumerate}
\item
If $X_i$ is not a bounded valence bushy tree, then $X'_i$ 
     is isometric to $X_{i'}$ for some $i'$ in the $G$-orbit 
$G(i)$ of $i$.
\item If $X_i$ is a bounded valence bushy tree, then so is $X'_i$.
\end{enumerate}
\end{theorem}
 As in the previous corollary, it is necessary to permit
$X_i'$ to be nonisometric to $X_i$.  Moreover, there may be
factors $X_i$ and $X_j$ of type (4) lying in the same $G$-orbit, but which 
are not even homothetic, so it is not sufficient to allow
rescaling of factors. 
\proof
We first assume that the action $G\acts I$ is transitive.  Pick $n\in I$.
Then the quasi-action $G_{n}\acts X_{n}$ is quasi-conjugate to an isometric
action $G_{n}\acts X_n'$, where $X_n'$ is isometric to $X_n$ unless $X_n$
is a bounded valence bushy tree, in which case $X_n'$ is a 
bounded valence bushy tree but not necessarily 
isometric to $X_n$; 
this follows from:

$\bullet$ \cite{hinkkanen,gabai,cassonjungreis,markovic}
when $X_n$ is $\H^2$.  Note that any quasi-action on $\H^2$ is
quasi-conjugate to an isometric action.

$\bullet$  \cite{sullivan,gromov,tukia,pansu,chow} when $X_n$ is a
rank $1$ symmetric space other than $\H^2$.  Note that Sullivan's theorem
implies  that any quasi-action on $\H^3$ is quasi-conjugate to an isometric 
action.  Also,  the proof given
in Chow's paper on the complex hyperbolic case covers arbitrary cobounded
quasi-actions, even though it is only stated for discrete cobounded quasi-actions.

$\bullet$ \cite{kl,lhabilitation} when $X_n$ is an irreducible symmetric space
or Euclidean building of rank at least $2$. 

$\bullet$ \cite{moshersageevwhyte} when $X_n$ is a bounded valence
bushy tree. 

By Theorem \ref{thminducedproperties}, the associated induced 
quasi-action of $G$ is 
quasi-conjugate
to the original quasi-action $G\acts \prod_{i\in I}\;X_i$ by a 
product quasi-isometry,
and we are done.

In the general case, for each orbit $G(i)\subset I$ of the 
action $G\acts I$, 
we have a well-defined associated quasi-action 
$G\acts \prod_{j\in G(i)}\; X_{j}$
for which the theorem has already been established, and we obtain the desired
isometric action $G\acts \prod_{i\in I}\;X_i'$  by taking products.  
\qed

\bigskip
\begin{corollary}
\label{corgroupcorollary}
Let $\{X_i\}_{i\in I}$ be as above, and suppose $G$ is a finitely
generated group quasi-isometric to the product $\prod_{i\in I}\;X_i$.  Then 
$G$  admits a discrete, cocompact, isometric
action on a product $\prod_{i\in I}\;X_i'$, where
for every $i$,  $X_i'$ is quasi-isometric to $X_i$,
and precisely one of the following holds:
\begin{enumerate}
\item 
$X_i$ is not a bounded valence bushy tree, and 
$X'_i$ is isometric to $X_{i'}$ for some $i'$ in the 
$G$-orbit $G(i)\subset I$ of $i$.
\item  Both  $X_i$ and $X_i'$ are   bounded valence bushy trees.
\end{enumerate}
\end{corollary}
\proof
Such a group $G$ admits a discrete, cobounded quasi-action on 
$\prod_{i\in I}\;X_i$.  Theorem \ref{thmquasiaction} furnishes the desired
 isometric action $G\acts \prod_i\; X_i'$.
\qed

\bigskip
{\em Remarks.}
\begin{itemize}
\item
Corollary \ref{corgroupcorollary} refines earlier results 
\cite{ahlin,gpsqiss,moshersageevwhyte}.
\item A proper Gromov hyperbolic space with cocompact isometry group 
is
of coarse type I unless it is quasi-isometric to $\R$ 
\cite[Sec. 3]{kkl}.
\item The classification of the four different types of spaces 
above is
quasi-isometry invariant, with one exception: a space of type (1) 
will
also be a space of type (4) iff it is a quasi-isometrically 
rigid rank $1$ symmetric space (i.e. a quaternionic hyperbolic 
space or the 
Cayley hyperbolic plane \cite{pansu}).  See Lemma \ref{lemtypeclassification}.
\item Two irreducible symmetric spaces are quasi-isometric iff they
are isometric, up to rescaling \cite{mostow,pansu,kl}.  Two Euclidean
buildings as in (2) above are quasi-isometric iff they are isometric 
up
to rescaling \cite{kl,lhabilitation}. 
\end{itemize}

\section{The construction of induced quasi-actions}

The construction of induced quasi-actions is a direct
imitation of one of the standard constructions of induced 
representations.   
We now review this for the convenience of the reader.

Let $H$ be a subgroup of some group $G$, and suppose
$\al:H\acts V$ is a linear representation.  Then we have an action
$H\acts G\times V$ where $(h,(g,v))=(gh^{-1},hv)$.  Let 
$E\defeq (G\times V)/H$ be the quotient.
 There is  a natural projection map  $\pi:E\ra G/H$ whose fibers are copies
of $V$;  this would be a vector bundle over the discrete space $G/H$
if $V$ were endowed with a 
topology.  The action $G\acts G\times V$ by left translation on the
first factor descends to $E$, and commutes with the projection
map $\pi$.  
Moreover, it preserves the linear structure on the fibers.
Hence there is a representation of $G$ on the vector
space of sections $\Ga(E)$, and this is the representation of $G$
induced by $\al$.

We use the terminology of \cite[Sec.\ 2]{gpsqiss}. 
(However, we replace {\em quasi-isometrically conjugate} by the shorter and 
more accurate term {\em quasi-conjugate}.)

We will work with generalized metrics taking values in $[0,+\infty]$. 
A {\em finite component} of a generalized metric space is an equivalence class
of points with pairwise finite distances.
Clearly, quasi-isometries respect finite components. 

Let $\{X_i\}_{i\in I}$ be a finite collection 
of unbounded metric spaces in the usual sense,
i.e.\ the metric on each $X_i$ takes only finite values.
On their product $\prod_{i\in I}\;X_i$
we consider the natural ($L^2$-)product metric. 
On their disjoint union $\sqcup_{i\in I}\;X_i$
we consider the generalized metric 
which induces the original metric on each component $X_i$ 
and gives distance $+\infty$ to any pair of points in different components. 

We observe that a quasi-isometry
$\prod_{i\in I}\;X_i\ra \prod_{i\in I}\;X'_i$
preserving the product structure 
gives rise to a quasi-isometry 
$\sqcup_{i\in I}\;X_i\ra \sqcup_{i\in I}\;X'_i$,
well-defined up to bounded error,
and vice versa.
Thus equivalence classes of quasi-actions 
$\al:G\acts\prod_{i\in I}\;X_i$
preserving the product structure 
correspond one-to-one to 
quasi-actions 
$\beta:G\acts\sqcup_{i\in I}\;X_i$.
In what follows we will prove  the disjoint union analog 
of Theorem \ref{thminducedproperties}. 
(The index of $H$ can be arbitrary from now on.) 

\begin{lemma}
Suppose that $Y$ is a generalized metric space 
and that $G\acts Y$ is a quasi-action such that 
$G$ acts transitively on the set of finite components of $Y$.
Let $Y_0$ be one of the finite components and $H$ its stabilizer in $G$.
Then the restricted action $H\acts Y_0$ 
determines the action $G\acts Y$ up to quasi-conjugacy.
\end{lemma}
\proof
If $G\acts Y'$ is another quasi-action,
$Y_0'$ is a finite component with stabilizer $H$, 
then any quasi-conjugacy between
$H\acts Y_0$ and $H \acts Y_0'$ extends in a straightforward way to a
quasi-conjugacy between $G \acts Y$  and $G \acts Y'$.
\qed

We will now show how to 
recover the $G$-quasi-action from the $H$-quasi-action
by quasifying the construction of induced actions 
as described above. 

\begin{definition}
An {\em $(L,A)$-coarse fibration} 
$(Y,{\mathcal F})$ consists of a (generalized) metric space $Y$ 
and a family ${\mathcal F}$ of subsets $F\subset Y$,
the {\em coarse fibers},  
with the following properties:
\begin{enumerate}
\item The union $\cup_{F\in{\mathcal F}}F$ of all fibers 
has Hausdorff distance $\leq A$ from $Y$.
\item For any two fibers $F_1,F_2\in{\mathcal F}$ holds 
\begin{equation*}
d_H(F_1,F_2) \leq L\cdot d(y_1,F_2)+A \qquad\forall\; y_1\in F_1.
\end{equation*}
\end{enumerate}
We also say that ${\mathcal F}$ is a coarse fibration of $Y$.
\end{definition}

Note that the coarse fibers are not required to be disjoint. 

It follows from part (2) of the definition 
that $d_H(F_1,F_2)<+\infty$ if and only if $F_1$ and $F_2$
meet the same finite component of $Y$.
We will equip the ``base space'' ${\mathcal F}$ with the Hausdorff metric.

\begin{lemma}
\label{qorbfib}
If $H\acts Y$ is an $(L,A)$-quasi-action 
then the collection of quasi-orbits $O_y:=H\cdot y$ 
forms an $(L,3A)$-coarse fibration of $Y$. 
\end{lemma}
\proof
For $h,h_1,h_2\in H$ and $y_1,y_2\in Y$ we have
\begin{equation*}
d\bigl(hy_1,(hh_1^{-1}h_2)y_2\bigr)) 
\leq 
d\bigl((hh_1^{-1})(h_1y_1),(hh_1^{-1})(h_2y_2)\bigr)) +2A
\leq
L\cdot d(h_1y_1,h_2y_2)+3A
\end{equation*}
and so 
\begin{equation*}
d(O_{y_1},O_{y_2})\leq L\cdot d(h_1y_1,O_{y_2})+3A .
\end{equation*}
\qed

Let $(Y,{\mathcal F})$ and $(Y',{\mathcal F}')$ be coarse fibrations. 
We say that a map $\phi:Y\ra Y'$ 
{\em quasi-respects} the coarse fibrations 
if the image of each fiber $F\in{\mathcal F}$ 
is uniformly Hausdorff close to a fiber $F'\in{\mathcal F}'$,
$d_H(\phi(F),F')\leq C$.
The map $\phi$ then induces a map $\bar\phi:{\mathcal F}\ra{\mathcal F}'$ 
which is well-defined up to bounded error $\leq2C$.
Observe that if $\phi$ is an $(L,A)$-quasi-isometry then
$\bar\phi$ is an $(L,A+2C)$-quasi-isometry. 

We say that a quasi-action $\rho:G\acts Y$ 
{\em quasi-respects} a coarse fibration ${\mathcal F}$ 
if all maps $\rho(g)$ quasi-respect ${\mathcal F}$ 
with uniformly bounded error.
The quasi-action $\rho$ then descends to a quasi-action 
$\bar\rho:G\acts{\mathcal F}$
which is unique up to equivalence 
(cf.\ \cite[Definition 2.3]{gpsqiss}). 

We apply these general remarks to the following situation
in order to obtain our main construction. 

Let $G$ be a group, $H<G$ a subgroup (of arbitrary index)
and $H\buildrel\al\over\acts X$ an $(L,A)$-quasi-action.
Let $Y=G\times X$ where $G$ is given the metric 
$d(g_1,g_2)=+\infty$ unless $g_1=g_2$.
That is, 
$Y$ consists of $|G|$ finite components each of which is a copy of $X$.  
The quasi-action $\al$ gives rise to a product quasi-action 
$H\buildrel\rho_H\over\acts Y$ via
\begin{equation*}
\rho_H\bigl(h,(g,x)\bigr) = (gh^{-1},hx) .
\end{equation*}
We denote by ${\mathcal F}_H$ the coarse fibration of $Y$
by $H$-quasi-orbits.
The isometric $G$-action given by 
\begin{equation*}
\tilde\rho_G\bigl(g',(g,x)\bigr) = (g'g,x) 
\end{equation*}
commutes with $\rho_H$.
As a consequence,
$\tilde\rho_G$ descends to an isometric action
\begin{equation}
\label{induced_union}
\hat\beta:=\bar\rho_G:G\acts{\mathcal F}_H .
\end{equation}
If $H=G$ then $\al$ is quasi-conjugate to $\hat\beta$
via the quasi-isometry $x\mapsto \rho_H(H)\cdot(e,x)$. 

In general, 
the finite components of ${\mathcal F}_H$ correspond 
to the left $H$-cosets in $G$. 
More precisely,
$gH$ corresponds to 
$\cup_{x\in X} \rho_H(H)\cdot(g,x)$,
that is,
to the union of $\rho_H$-quasi-orbits contained in $gH\times X$. 
$H$ stabilizes the finite component
$\cup_{x\in X} \rho_H(H)\cdot (e,x)$.
The action of $H$ on this component 
is quasi-conjugate to $\al$.

As remarked in the beginning of this section,
$\hat\beta$ is the unique $G$-quasi-action up to quasi-conjugacy 
such that $G$ acts transitively on finite components
and such that $H$ is the stabilizer of a finite component 
and the restricted $H$-quasi-action is quasi-isometrically conjugate to $\al$. 

Passing back from disjoint unions to products
we obtain Theorem \ref{thminducedproperties}.

\section{Quasi-isometries and the classification into types (1)-(4)}

 We now prove:

\begin{lemma}
\label{lemtypeclassification}
Suppose $Y$ and $Y'$ are  spaces of one of types (1)-(4) 
as in Theorem \ref{thmquasiaction}.  If $Y$ is quasi-isometric to 
$Y'$,
then they have the same type, unless one is a 
quasi-isometrically rigid rank $1$ symmetric space,
and the other is of type (4).
\end{lemma}
\proof
First suppose one of the spaces is not Gromov hyperbolic.  Since Gromov hyperbolicity
is quasi-isometry invariant, both spaces must be higher rank space of either of type (1) or (2).
But by \cite{kl}, two irreducible symmetric spaces or Euclidean buildings 
of rank at least two are quasi-isometric iff they are homothetic.  Thus in this
case they must have the same type.

Now assume both spaces are Gromov hyperbolic.  Then $\geo Y$ and $\geo Y'$
are homeomorphic. 

 If $Y$ is a bounded valence bushy tree, then it is well-known
that $Y$ is quasi-isometric to a trivalent tree, and $\geo Y$
is homeomorphic to a Cantor set.  Therefore $Y$ cannot be quasi-isometric
to a space of type (1), since the boundary of a Gromov hyperbolic symmetric
space is a sphere.  Also, the quasi-isometry group of a 
trivalent tree $T$ has an induced action on the space of triples in $\geo T$ which is
not proper, and hence
it cannot be quasi-isometric to a space of type (4).

If $Y$ is a hyperbolic or complex hyperbolic space, then the induced
action of $\qi(X)$ on the space of triples in $\geo X$ is not proper,
and hence $Y$ cannot be quasi-isometric to a space of type (4).

The lemma follows.
\qed

\bibliography{induced}
\bibliographystyle{alpha}

\end{document}